\documentclass[12pt]{article}
\usepackage{amsmath,amssymb}
\usepackage{graphicx}

\newcounter{tally}

\newtheorem {lemma}{Lemma}
\newtheorem {corollary}{Corollary}

\newtheorem {thm}{Theorem}

\def \a   {\alpha}
\def \b   {\beta}
\def \g   {\gamma}
\def \d   {\delta}

\def \D   {{\Delta}}
\def \eps {\varepsilon}

\def\E{{\mathbb{E\,}}}

\def\P{{\mathbb{P}}}
\def \R {{\mathbb{R}}}
\newcommand{\Z}{\mathbb{Z}}
\def\F{{\cal{F}}}
\def\|{\,|\,}

\def\limt{\lim_{t\to\infty}}
\def\limt0{\lim_{t\to 0}}

\def\bn{\begin{eqnarray*}}
\def\en{\end{eqnarray*}}
\def\bnn{\begin{eqnarray}}
\def\enn{\end{eqnarray}}

\newcommand{\proof}{{\sf Proof.\ }}
\newcommand{\BBox}{\rule{6pt}{6pt}}
\newcommand\Cox{$\hfill \BBox$ \vskip 5mm}

\def \Rec   {{\sf \bf Rec}}
\def \Trans   {{\sf \bf Trans}}

\title{Urn-related random walk with drift $\rho \,x^\a\,/\,t^\b$}

\author{Mikhail Menshikov\footnote{Department of Mathematics, University of Durham, DH1~3LE, U.K.\newline
E-mail:~Mikhail.Menshikov@durham.ac.uk@durham.ac.uk} and Stanislav
Volkov\footnote{Department of Mathematics, University of Bristol,
BS8~1TW, U.K.
\newline
E-mail:~S.Volkov@bristol.ac.uk}}
\date{9 November  2007}

\begin {document}

\maketitle
\begin {abstract}
We study a one-dimensional random walk whose expected drift
depends both on time and the position of a particle. We establish
a non-trivial phase transition for the recurrence vs.\ transience
of the walk, and show some interesting applications to Friedman's
urn, as well as showing the connection with Lamperti's walk with
asymptotically zero drift.
\end {abstract}

\noindent {\bf Keywords:} random walks, urn models, martingales.

\noindent {\bf Subject classification:} primary 60G20; secondary
60K35.

\section{Introduction}\label{Intro}
%
%
%
Consider the following  stochastic processes $X_t$ which may
loosely be described as a random walk on $\R_+$ (or in more
generality on $\R$) with the asymptotic drift given by
 \bn
 \mu_t:=\E (X_{t+1}-X_t\| X_t=x)\sim \rho \frac
 {|x|^\alpha}{t^\beta},
 \en
where $\rho$, $\a$ and $\b$ are some fixed constants, and the
exact meaning of ``$\sim$'' will be made precise later. In this
paper we establish when this process is recurrent or transient, by
finding the whole line of phase transitions in terms of $(\a,\b)$.
We also analyze some critical cases, when the value of $\rho$
becomes important as well. Note that because of symmetry, it is
sufficient to consider only these processes on $\R_+$, and from
now on we will assume that $X_t\ge 0$ a.s.\ for all values of $t$.

The original motivation of this paper is based on an open problem
related to Friedman urns, posed in Freedman~\cite{Fr}. In certain
regimes of these urns, to the best of our knowledge, it is still
unknown whether the number of balls of different colors can
overtake each other infinitely many times with a positive
probability. We will not describe this problem in more details
here, rather we refer the reader directly to
Section~\ref{sec_urn}.

Incidentally, the class of stochastic processes we are considering
covers simultaneously not only the Friedman urn, but also the walk
with an asymptotically zero drift, first probably studied by
Lamperti, see~\cite{Lam1} and~\cite{Lam2}. His one-dimensional
walks with drift depending only on the position of the particle
naturally arise when proving recurrence of the simple random walk
on $\Z^1$ and $\Z^2$ and transience on $\Z^d$, $d\ge 3$. They can
be used of course for answering the question of recurrence for a
much wider class of models, notably those involving polling
systems, for example, see~\cite{AIM}  and~\cite{MZ}. It will be
not surprising if the model we are considering also covers some
other probabilistic models, of which we are unaware at the moment.

%

\begin{figure}[hbt]
     \centering
     \unitlength=0.35mm
     \begin{picture}(400.00,200.00)
          \put(0,0){\line(1,0){400.00}}
           \put(200,0){\dashbox{5}(0,200)}
          \put(200,0){\line(1,1){150,150}}
          \put(100,0){\line(2,1){200,200}}
          \put(200,0){\vector(1,0){200}}
          \put(200,200){\vector(0,1){10}}
          \put(300,100){\circle*{5.00}}
          \put(200,-8){\makebox(0,0)[cc]{$0$}}
          \put(100,-8){\makebox(0,0)[cc]{$-1$}}
          \put(300,-8){\makebox(0,0)[cc]{$1$}}
          \put(190,55){\makebox(0,0)[cc]{$\frac 12$}}
          \put(190,100){\makebox(0,0)[cc]{$1$}}
          \put(390,-8){\makebox(0,0)[cc]{$\alpha$}}
          \put(190,200){\makebox(0,0)[cc]{$\beta$}}
          \put(190,150){\makebox(0,0)[cc]{recurrence}}
          \put(200,30){\makebox(0,0)[cc]{transience}}
          \put(300,30){\makebox(0,0)[cc]{prohibited area}}
          \put(320,100){\makebox(0,0)[cc]{$(1,1)$}}
     \end{picture}
     \caption{Diagram for $(\alpha,\beta)$.
     \label{SelfdualFig}}
\end{figure}
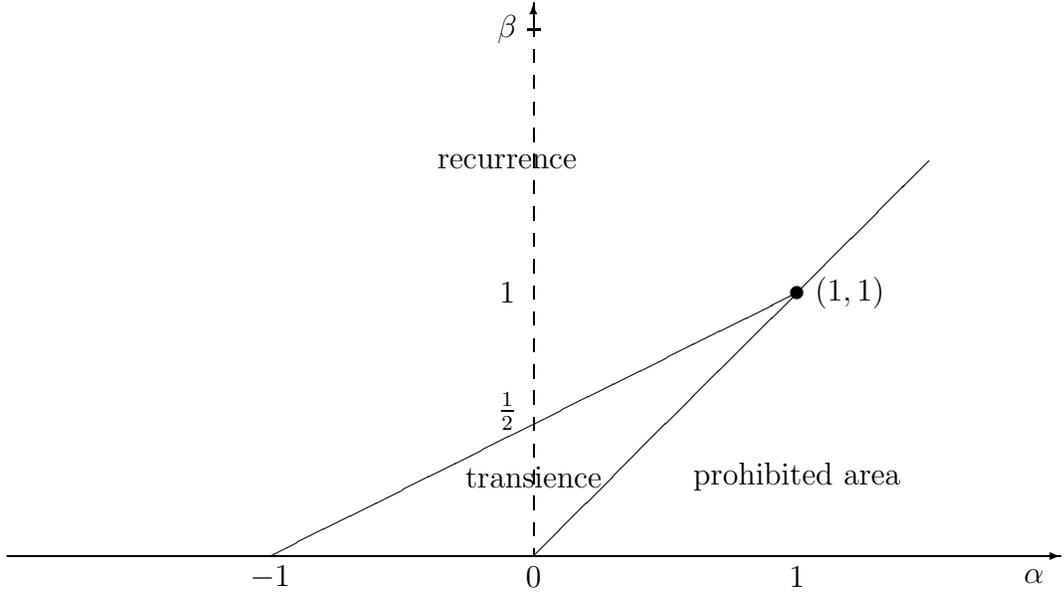

In our paper we study the random walk whose drift depends {\sl
both on time and the position} of a particle. Throughout the paper
we assume that
$$
 (\a,\b)\in \Upsilon=\{(\a,\b):\ \b > \a \text{ and } \b\ge 0 \}
$$
to avoid the situations when the drift becomes unbounded and the
borderline cases (the only exception will be $\a=\b=1)$. We will
show that under some regularity conditions, the walk is transient
when $(\a,\b)$ lie in the following area
$$
 \Trans=\left\{(\a,\b):\ 0\le\b< 1, \ 2\b-1<\a <\b \right\}\subset \Upsilon
$$
and recurrent for $(\a,\b)$ in
  \bn
 \Rec&=&\Upsilon \setminus \overline{\Trans}=\{(\a,\b):\ \b\ge 0, \  \a<\min(\b,2\b-1)\}
 \en
where $\overline{\Trans}$ denotes the closure of the set $\Trans$.
In the special critical case $\a=\b=1$ we show that the walk is
transient for $\rho>1/2$ and recurrent for $\rho<1/2$. An example
of such a walk with $\a=\b=1$ is the process on $\Z_+$ with the
following jump distribution:
 \bn
 \P(X_{t+1}=n\pm 1\| X_t=n)=\frac 12
\pm \frac{\rho n}{2t}
  \en
This walk is analyzed in Section~\ref{sec_urn}.

Throughout the paper we will need the following hypothesis. Let
$X_t$ be a stochastic process on $\R_+$ with jumps
$D_t=X_t-X_{t-1}$ and let $\F_t=\sigma(X_0,X_1,\dots,X_t)$. Let
$a$ be some positive constant.
\\ \vskip 1mm \noindent
 \quad {\bf (H1)} \qquad Uniform boundedness of jumps. \\ There is a constant $B_1>0$ such that $|D_t|\le B_1$ for all $t\in\R_+$ a.s.
\\ \vskip 1mm \noindent
 \quad {\bf (H2)} \qquad Uniform non-degeneracy on $[a,\infty)$.\\ There is a constant $B_2>0$ such that whenever $X_{t-1}\ge a$,
  $\E (D_t^2\|\F_{t-1})\ge  B_2$ for all $t\in\R_+$.
 a.s.
\\ \vskip 1mm \noindent
 \quad {\bf (H3)} \qquad
 Uniform boundedness of time to leave [0,a].\\
 The number of steps required for $X_t$ to exit the interval $[0,a]$ starting from any point inside
 this interval is uniformly stochastically bounded above by some
 independent random variable $W\ge 0$ with a finite mean $\mu=\E W<\infty$, i.e., for all
 $s\in\R_+$, when $X_s\le a$,
 \bn
 \forall x\ge 0 \ \ \P(\eta(s)\ge x\|\F_s)\le \P(W\ge
 x),
 \mbox{where }\eta(s)=\inf\{t\ge s:\ X_t>a\}.
 \en
\\ \vskip 1mm \noindent
 The rest of the paper is organized as follows. In
Section~\ref{sec_tech} we prove some technical lemmas. In
Section~\ref{sec_tr} we formulate the exact statement about the
transience of the process $X_t$ and prove it while in
Section~\ref{sec_rec} we do the same for recurrence. We also study
some borderline cases in Section~\ref{sec_spec}, and present an
open problem in Section~\ref{subsec_open}. Finally, we apply our
results to generalized P\'olya and Friedman urns in
Section~\ref{sec_urn}.

\section{Technical facts}\label{sec_tech}

First, we will need the following important claim about the law of
iterated logarithms for martingales.
\begin{lemma}[Proposition (2.7) in Freedman (1975)] \label{lem_mart}
Suppose that $S_n$ is a martingale adapted to filtration $\F_n$
and $\D_n=S_n-S_{n-1}$ are its differences. Let $T_n=\sum_{i=1}^n
\E(\D_i^2\|\F_{i-1})$,  $\sigma_b=\inf\{n:\ T_n>b\}$, and
 \bn
 L(b)=ess\sup_{\omega} \sup_{n\le \sigma_b(\omega)}|\D_n(\omega)|.
 \en
Suppose $L(b)=o(b/\log\log b)^{1/2}$ as $b\to\infty$. Then
 \bn
 \limsup_{n\to\infty} \frac{S_n}{\sqrt{2T_n\log\log T_n}}=1
  \ \text{ a.s.\ on } \{T_n\to\infty\}.
 \en
\end{lemma}

\begin{lemma}\label{lem_lil}
Let $X_t$, $t=1,2,\dots$ be a sequence of random variables adapted
to filtration $\F_t$ with differences $D_t=X_{t}-X_{t-1}$
satisfying (H1),  (H2), and (H3) for some $a>0$. Suppose that on
the event $\{X_t\ge a\}$
 \bn
 &\E (D_{t+1}\|\F_{t})\ge 0 \mbox{ a.s.}
 \en
Then for any $A>0$
$$
\P(\exists t:\ X_t>A \sqrt t)=1.
$$
\end{lemma}
Note that in Lemma~\ref{lem_lil} we prove a weaker result than in
the original Lemma~\ref{lem_mart}. This result, one hand hand,
will suffice for our purposes, while on the other hand it does not
require a sequence  $X_t$ to be an exact martingale, rather than
just a submartingale, and only on $[a,\infty)$. \vskip 1mm

 \proof
First, we are going essentially to ``freeze'' the process $X_t$
whenever it enters the interval $[0,a]$, where it is not a
submartingale, until the moment when $X_t$ exits from this
interval.
Define the function $s(t):\Z_+\to\Z_+$ such that $s(0)=0$ and for
$t\ge 0$
 \bn
 s(t+1)=\left\{\begin{array}{ll}
 t+1, & \mbox{ if $X_t>a$ or $X_{t+1}>a$}\\
 s(t), & \mbox{ otherwise.}
 \end{array}\right.
 \en
Let $\tilde X_t=X_{s(t)}$. Then $\tilde X_n$ is a submartingale
satisfying (H1), perhaps with a new constant $\tilde B_1=B_1+a$.
Indeed, when $X_t>a$, $\tilde X_t=X_t$ and $\tilde
X_{t+1}=X_{t+1}$, so $\E(\tilde X_{t+1}-\tilde
X_t|\F_t)=\E(D_{t+1}|\F_t)\ge 0$. When $X_t<a$ (and so is $\tilde
X_t<a$), either $X_{t+1}<a$ and then $s(t+1)=s(t)$ implying
$\tilde X_{t+1}=\tilde X_t$, or $X_{t+1}\ge a$ in which case
$\tilde X_{t+1}=X_{t+1}\ge a> \tilde X_t$.

Let
 \bn
 \tilde D_n&=&\tilde X_n-\tilde X_{n-1}\\
 Z_n&=&\E(\tilde D_n\|\F_{n-1})\ge 0\\
 S_n&=&X_n-Z_1-Z_2-\dots-Z_{n}.
 \en
Then
$$
\E(S_{n}-S_{n-1}\|\F_{n-1})=\E(X_{n}-X_{n-1}-Z_n\|\F_{n-1})=0
$$
whence $S_n$ is a martingale with differences
$\D_n:=S_{n}-S_{n-1}=\tilde D_{n}-Z_n$. Note that
 \bnn\label{eq_DD}
 \E(\D_n^2\|\F_{n-1})=\E((S_{n}-S_{n-1}-Z_{n-1})^2\|\F_{n-1})
  =  \E(\tilde D_n^2\|\F_{n-1})-Z_n^2.
 \enn

Let $\eta_0=0$ and for $k=1,2,\dots$
 \bn
 \zeta_k&=&\inf\{t\ge \eta_{k-1}:\ X_t\le a\},\\
  \eta_k&=&\inf\{t\ge \zeta_{k} :\ X_t> a\}
 \en
be the consecutive times of entry in and exit from $[0,a]$.  Then
$\tilde W_k:=\eta_k-\zeta_k$ are stochastically bounded by i.i.d.\
random variables $W_1,W_2,\dots$ with the distribution of $W$.
Therefore
$$
 \limsup_{m\to\infty} \frac{\sum_{i=1}^m\tilde W_i}{m}\le \mu\mbox{ a.s.}
$$
and consequently the number
 \bn
 I_n&=&\{t\in\{0,1,\dots,n-1\}: t\not\in[\zeta_k,\eta_{k})\mbox{ for some }k\}
 \\
 &=&\{t\in\{0,1,\dots,n-1\}: X_t>a\}
 \en
of those times which do not belong to some ``frozen'' interval
$[\zeta_k,\eta_k)$  satisfies a.s.
 \bnn\label{eq_In}
 |I_n| \ge \frac n{2\mu}
 \enn
for $n$ sufficiently large.

Next, since $\tilde D_n$'s are bounded, we have $|\D_n|\le |\tilde
D_n|+|\E(\tilde D_n\|\F_{n-1})|\le 2(B_1+a)$. Therefore, $L(b)\le
2(B_1+a)$ and the conditions of Lemma~\ref{lem_mart} are met.
First, suppose that
$$
T_n=\sum_{i=1}^n \E(\D_i^2\|\F_{i-1})\to\infty,
$$
then
 \bn
 \limsup_{n\to\infty} \frac{S_n}{\sqrt{2T_n\log\log T_n}}=1 \ \mbox{ a.s.}
 \en
Therefore, for infinitely many $n$'s we would have
 \bn
 S_n\ge \sqrt{T_n\log\log T_n}.
 \en
Using (\ref{eq_DD}), this results in
 \bnn\label{eq_xz}
 X_n&=&\sum_{i=1}^n Z_i+ S_n\ge \sum_{i=1}^n Z_i+ \sqrt{ \sum_{i=1}^n(\E(\tilde D_i^2\|\F_{i-1})-Z_i^2)\log\log T_n}
 \nonumber \\
 &\ge&
 \sum_{i\in 1+I_n}^n Z_i+ \sqrt{ \sum_{i\in 1+I_n}(\E(D_i^2\|\F_{i-1})-Z_i^2)\log\log
 T_n}
 \enn
since $i-1\in I_n$ implies $X_{i-1}>a$ and consequently $\tilde
D_{i}=D_{i}$ (note that each term in the sums above is
non-negative). Let $0\le k\le |I_n|$ be the number of those
$Z_i$'s, $i\in I_n$ such that $Z_i<\sqrt{B_2/2}$.
Then~(\ref{eq_xz}) together with $\E(D_i^2\|\F_{i-1})\ge B_2$
yields
 \bn
  X_n  \ge (|I_n|-k)\sqrt{\frac{B_2}2}+ \sqrt{\frac {kB_2\log\log T_n}2}
  \ge \sqrt{\frac {nB_2\log\log  T_n}{2\mu}}
 \en
for $n$ large enough, taking into account the fact that $T_n\le
B_1^2 n$ and inequality~(\ref{eq_In}). This implies
 the statement
of Lemma~\ref{lem_lil}, since we assumed $T_n\to\infty$.

On the other hand, on the complementary event $\sum_{i=1}^{\infty}
\E(\D_i^2\|\F_{i-1})<\infty$, by e.g.\ Theorem in Chapter 12 in
Williams (1991) $S_n$ converges a.s.\ to a finite quantity
$S_\infty$, and we obviously must also have $\E(\tilde
D_n^2\|\F_{n-1})-Z_n^2\to 0$ yielding
 \bn
\liminf_{i\to\infty,\ i:\ X_{i-1}>a} Z_i\ge \sqrt{B_2} .
 \en
Combining
this with~(\ref{eq_In}), we obtain
 \bn
\liminf_{n\to\infty} \frac{X_n}n=\liminf_{n\to\infty} \frac{S_n
+Z_1+Z_2+\dots+Z_n}n\ge \frac{\sqrt{B_2}}{2\mu}
 \en
which is even a stronger statement than we need to prove.
 \Cox

\begin{lemma}\label{lem_snos}
Fix $a>0$, $c>0$, $\g\in(0,1)$, and consider a Markov process
$X_t$, $t=0,1,2,\dots$ on $\R_+$ with jumps $D_t=X_{t}-X_{t-1}$,
for which the hypotheses (H1) and (H2) hold. Suppose that for some
large $n>0$ the process starts at
 $X_0\in(a,\g n ]$,
and that on the event $\{a\le X_t\le n\}$
$$
\E(D_{t}\| \F_{t-1})\le \frac cn.
$$
Let
$$
\tau=\inf\{t: \ X_t<a\text{ or } X_t> n\}.
$$
be the time to exit $[a,n]$. Then
\begin{list} {(\roman{tally})} {\usecounter{tally}}
 \item $\tau<\infty$ a.s.;
 \item $\P(X_{\tau}<a)\ge \nu=\nu(\g,c,B_2)>0$ uniformly in $n$.
\end{list}
\end{lemma}
\proof First, let us show that the process $X_t$ must exit $[a,n]$
in a finite time. Since $|D_t|\le B_1$, by Markov inequality for
non-negative random variables for any $\eps>0$ we have
 \bnn\label{eq_bd}
\P( B_1^2 - D_t^2 \ge  (1-\eps^2)B_1^2 \| \F_{t-1}) \le
\frac{\E(B_1^2 - D_t^2\| \F_{t-1})}{(1-\eps^2)B_1^2}\le
(1-\eps^2)^{-1}\left(1-\frac{B_2}{B_1^2}\right)
 \enn
Hence for a sufficiently small $\eps>0$ the RHS of (\ref{eq_bd})
can be made smaller than 1, whence there is a $\d>0$ such that
$$
\P( D_t^2 \ge  (\eps B_1)^2 \| \F_{t-1})> 2\d.
$$
In turn, this implies that at least one of the probabilities $\P(
D_t \ge  \eps B_1 \| \F_{t-1})$ or $\P( D_t \le  -\eps B_1 \|
\F_{t-1})$ is larger than $\d$. Hence from any starting point the
walk can exit $[a,n]$ in at most $n/(\eps B_1)$ steps with
probability at least $\d^{n/(\eps B_1)}$, yielding that $\eps B_1
\tau/n$ is stochastically bounded by a geometric random variable
with parameter $\d^{n/(\eps B_1)}$, which is not only finite but
also has all finite moments.

To prove the second claim of the lemma, first we establish the
following elementary inequality. Fix a $k\ge 1$ and consider the
function $g(x)=(1-x)^k-1+kx-k(k-1)x^2/4$. Since $g(0)=0$,
$g'(0)=0$, and $g''(x)=k(k-1)((1-x)^{k-2}-1/2)\ge 0$ for $|x|\le
1/(2k)$, we have $g(x)\ge 0$ on this interval. Consequently,
 \bnn\label{eq_elin}
 (1-x)^k -1 \ge -kx+\frac {k(k-1) x^2}4 \mbox{ for }x\in \left[-\frac1{2k},\frac1{2k}\right].
 \enn
Now let $Z_t=2n-X_t$ and $Y_t=Z_t^k$ for some $k\ge 1$ to be
chosen later. Suppose that $n>2kB_1$. Then, on the event $\{X_t\in
[a,n]\}$ we have $Z_t\in [n,2n]$ yielding $|D_{t+1}/Z_t|\le
B_1/n\le 1/(2k)$ and thus by (\ref{eq_elin}) we have
 \bn
 \E(Y_{t+1}-Y_{t}\|\F_{t})&=& Y_t
 \E\left(\left(1-\frac{D_{t+1}}{Z_t}\right)^k-1\|\F_{t}\right)
 \\&\ge&
 kY_t\left[
 -\frac{\E(D_{t+1}\|\F_{t})}{Z_t}+\frac{(k-1)\E(D_{t+1}^2\|\F_{t})}{4Z_t^2} \right]
\\
 &\ge &
kY_t\left[
 -\frac{c}{n Z_t}+\frac{(k-1)B_2}{4Z_t^2} \right]  \ge kY_t\left[
 -\frac{c}{n^2}+\frac{(k-1)B_2}{16n^2} \right]>0,
 \en
once $k>1+16c/B_2$.

Hence $Y_{t\wedge \tau}$ is a non-negative submartingale. By the
optional stopping theorem,
$$
\E(Y_\tau)\ge  Y_0\ge[(2-\g)n]^k.
$$
On the other hand,
$$
\E(Y_\tau)= \E(Y_\tau;\ X_\tau<a) + \E(Y_\tau;\ X_\tau>n)  \le
(2n)^k \P(X_{\tau}<a)+n^k (1-\P(X_{\tau}<a)),
$$
yielding
$$
\P(X_{\tau}<a) \ge \frac{(2-\g)^k-1}{2^k-1}=:\nu>0.
$$
\Cox

\begin{lemma}\label{lem_mart_clt}
Suppose that $X_t$, $t=0,1,\dots$ is a submartingale satisfying
(H1).  Then for any $x>0$
 \bn
  \P\left(\inf_{0\le t\le h x^2} X_t < X_0-bx\right)\le
  c(h,b,B_1)=\frac{4hB_1^2}{b^2}.
 \en
\end{lemma}
\proof Let $Z_n=\E(X_{n+1}-X_n\|\F_n)\ge 0$. Then
$$
S_t=X_0-(X_t-Z_1-Z_2-\dots-Z_t)=(X_0-X_t)+Z_1+\dots+Z_t\ge X_0-X_t
$$
is a square-integrable martingale with $S_0=0$, since $|S_n|\le
|X_0|+2nB_1$. Moreover, since
 \bn
\E\left((S_t-S_{t-1})^2\|\F_{t-1}\right) &=&
 \E\left((X_t-X_{t-1}-Z_t)^2\|\F_{t-1}\right)\\
 &=&\E\left((X_t-X_{t-1})^2\|\F_{t-1}\right)-Z_t^2\le B_1^2
 \en
we have
$$
A_n:=\sum_{t=1}^n \E\left((S_t-S_{t-1})^2\|\F_{t-1}\right) \le n
B_1^2.
$$
By Doob's maximum $\L^2$ inequality (see Durrett, pp.~254--255),
$$
\E\left(\sup_{0\le m \le n} |S_m|^2\right)\le 4 \E S_n^2=4A_n\le 4
n B_1^2.
$$
Consequently, by Chebyshev's inequality
 \bn
 \P\left(\inf_{0\le t\le h x^2} X_t<X_0-bx\right)
 &=&\P\left(\sup_{0\le t\le h x^2} X_0-X_t>bx\right)
  \\ &\le&
 \P\left(\sup_{0\le t \le hx^2} |S_t|> bx\right)
 < \frac{4(h x^2)B_1^2 }{b^2x^2}=\frac{4hB_1^2}{b^2}.
 \en
\Cox

\section{Transience}\label{sec_tr}
\begin{thm}\label{t1}
Consider a Markov process $X_t$, $t=0,1,2,\dots$ on $\R_+$ with
increments $D_t=X_{t}-X_{t-1}$ which satisfies (H1), (H2), and
(H3) for some $a>0$. Suppose that for $t$ sufficiently large on
the event $\{X_t\ge a\}$ we have either
\begin{list} {(\roman{tally})} {\usecounter{tally}}
 \item for some $\rho>1/2$
 \bn
 \E(D_{t+1}\| \F_{t})\ge \frac{\rho X_t}t,
 \en
 or
 \item for some $\rho>0$ and $(\a,\b)\in \Trans$
 \bn
 \E(D_{t+1}\| \F_{t})\ge \frac{\rho X_t^{\a}}{t^\beta}.
 \en
\end{list}
Then $X_t$ is transient in the sense that for any starting point
$X_0=x$  we have
$$
\P(\lim_{t\to\infty} X_t=\infty)=1.
$$
\end{thm}
\proof Consider  $Y_t=t/X_t^2$. Then
 \bn
Y_{t+1}-Y_t&=&\frac{t+1}{(X_t+D_{t+1})^2}-\frac t{X_t^2}  =
 \frac{t+1}{X_t^2}\left[
 \frac{1}{(1+D_{t+1}/X_t)^2}- \frac 1{1+1/t}
 \right]\\
 &\le&
\frac{t+1}{X_t^2}\left[
 \frac 1t-2\frac{D_{t+1}}{X_t}+3\frac{D_{t+1}^2}{X_t^2}+O\left(\frac{D_{t+1}}{X_t}\right)^3
 \right]\\
%
  \en
yielding
 \bnn\label{eq_supe}
\E(Y_{t+1}-Y_t\|\F_t)\le
 \frac
 {1+1/t}{X_t^2} Q_t
  \enn
where
$$
 Q_t=1-2\rho\frac{t^{1-\b}}{X_t^{1-\a}}+3B_1^2\frac{t}{X_t^2}+O(X_t^{-3}).
$$
Consider two cases:
\begin{list} {(\roman{tally})} {\usecounter{tally}}
 \item $\a=\b=1$, then $Q_t=1-2\rho+3B_1^2\frac{t}{X_t^2}+O(X_t^{-3})$;
 \item $(\a,\b)\in\Trans$.
\end{list}
In the first case, $Q_t$ and hence (\ref{eq_supe}) are negative as
long as $Y_t=t/X_t^2\le r$ for some  positive constant
$r<(2\rho-1)/3B_1^2$. (Note that this would imply $X_t\ge
\sqrt{t/r}\ge a$ for large enough $t$). Fix an arbitrary small
$\eps>0$ and suppose that for some time $s$ we have
$Y_s=s/X_s^2\le \eps r$. Let
$$
\tau=\tau(s)=\inf\{t>s:\ Y_t\ge r\}.
$$
Then $Y_{t\wedge \tau}$ is a non-negative supermartingale, hence
it a.s.\ converges to some random limit
$Y_\infty=\lim_{t\to\infty} Y_t$. By Fatou lemma, $\E Y_\infty\le
Y_s\le \eps r$. On the other hand,
$$
 \E Y_\infty=\E (Y_\infty;\ \tau<\infty)+\E (Y_\infty;\ \tau=\infty)\ge r\P(\tau<\infty)
$$
hence $\P(\tau<\infty)\le \eps$.

Finally, to show that for any $\eps>0$ with probability $1$ there
is an $s$ such that $s/X_s^2\le \eps r$ we apply
Lemma~\ref{lem_lil}. Consequently, $\P(\tau(s)=\infty\mbox{ for
some }s)=1$ yielding $\limsup_{t\to\infty} t/X_t^2\le r$ a.s., and
thus $\P(X_t\to\infty)=1$.

\vskip 5mm

Now consider case (ii) and observe that $0\le \b<1$ and $1-\a> 0$.
Suppose that
$$
 X_t^{2-2\d}\le t\le X_t^2 \text{ for some }
 \d\in\left(0,\frac{1+\a-2\b}{2(1-\b)}\right).
$$
Then
 \bn
 Q_t&=&1-2\rho\frac{t^{1-\b}}{X_t^{1-\a}}\left(1-\frac{3B_1^2}{2\rho}\frac{t^\b}{X_t^{1+\a}}\right)+O(X_t^{-3}).
 \en
Since $t\le X_t^2$, and $2\b<\a+1$,
 \bn
 \frac{t^{\b}}{X_t^{1+\a}}\le  \frac{X_t^{2\b}}{X_t^{1+\a}}=\frac
 1{X_t^{1+\a-2\b}}=o(1),
 \en
therefore
 \bn
 Q(t) &\le& 1-2\rho\frac{X^{2(1-\b)(1-\d)}}{X_t^{1-\a}}\left(1-o(1)\right)
  +O(X_t^{-3})= 1-2 \rho X^{2(1-\b)(1-\d)-(1-\a)} (1-o(1))<0
 \en
since $2(1-\b)(1-\d)-(1-\a)>0$ due to the choice of $\d$.
Therefore, on the event $\{ X_t^{2-2\d}\le t\le X_t^2\}$, $Y_t$ is
a supermartingale by inequality (\ref{eq_supe}).

Define the following areas:
 \bn
  M&=&\{s,x\ge 0: \ x^{2-\d}>s\},\\
  R&=&\{s,x\ge 0: \ x^{2-2\d}>s\},\\
  L&=&\{s,x\ge 0: \ x^{2}<s\}.
 \en
By Lemma~\ref{lem_lil}, there will be infinitely many times $s$
for which $s\le X_s^2/2$, so that $(s,X_s)\not\in L$. Fix such an
$s$ and let
$$
\tau=\tau(s)=\inf\{t>s:\ (t,X_t)\in L\cup M\}.
$$
Then $Y^{(*)}:=Y_{t\wedge\tau(s)}$ is a bounded supermartingale
which a.s.\ converges to $Y^{(*)}_\infty$; we have $\E
Y^{(*)}_\infty\le 1/2$ and as before obtain that on the event
$\{\tau<\infty\}$, $\P(Y_{\tau}\in L)\le 1/2$ independently of
$s$. Therefore, either $\tau(s)=\infty$ for some $s$ which implies
transience immediately, or by Borel-Cantelli lemma there will be
infinitely many times $s$ for which $(s,X_s)\in M$. From now
assume that the latter is the case.

Consider the sequence of stopping times when $(t,X_t)$ crosses the
curve $t=X_t^{2-\d}$, then reaches either area $L$ or area $R$
before crossing this curve again. Rigorously, suppose that for
some $t=\sigma_0$ we have $(t,X_t)\in M$ and it has just entered
area $M$. Set
$$
\eta_0=\inf\{t>\sigma_0:\ (t,X_t)\in L\cup R\}.
$$
Then for $k\ge 0$  let
 \bn
  \sigma_{k+1}&=&\left\{\begin{array}{ll}
   \inf\{t>\eta_k:\ (t,Y_t)\in M\},    & \mbox{if }(\eta_k,X_{\eta_k})\in L\\
   \inf\{t>\eta_k:\ (t,Y_t)\notin M\}, & \mbox{if }(\eta_k,X_{\eta_k})\in   R.
   \end{array}\right.
 \\
  \eta_{k+1}&=&\inf\{t>\sigma_{k+1}:\ (t,Y_t)\in L\cup R\}.
 \en
Thus we have
$$
 \sigma_0<\eta_0<\sigma_1<\eta_1<\sigma_2<\eta_2\dots
$$
Of course, it could happen that one of these stopping times is
infinity and hence all the remaining ones equal infinity as well;
however this would imply that $(t,X_t)\notin L$ for all large $t$,
which in turn implies transience (recall that we have assumed that
we visit the area $M$ infinitely often). Therefore, let us assume
from now on that all $\eta_k$'s and $\sigma_k$'s are finite.

For $t\ge \sigma_k$, $k\ge 0$, consider a supermartingale
$Y^{(k)}=Y_{t\wedge \eta_k}$. Since the jumps of $X_t$ are
bounded, $X_{\sigma_k}=\sigma_k^{1/(2-\d)}+O(1)$ and $\E
Y_{\eta_k}\le Y_{\sigma_k}=X_{\sigma_k}^{-\d}(1+o(1))$ and as
before, we obtain that
 \bnn\label{eq_p_l}
\P\left((\eta_k,X_{\eta_k})\in L \| \F_{\sigma_k}\right)\le
X_{\sigma_k}^{-\d} (1+o(1))=\frac{1+o(1)}{\sigma_k^{\d/(2-\d)}}.
 \enn
On the other hand, starting at $(\sigma_k,X_{\sigma_k})$ it takes
a lot of time for $(t,X_t)$ to reach $L$, and also if
$(\eta_k,X_{\eta_k})\in R$ it takes a lot of time to exit $M$,
since the walk has to go against the drift. More precisely,
 \bnn\label{eq_sigma_large}
 \sigma_{k+1}-\sigma_k
  \ge (\eta_k-\sigma_k)1_{(\eta_k,X_{\eta_k})\in L} +
   (\sigma_{k+1}-\eta_k)1_{(\eta_k,X_{\eta_k})\in R}.
 \enn
Set $x:=X_{\sigma_k}$,
$$
h=\frac 1{2(2B_1+1)^2}<\frac 1{8B_1^2},
$$
and observe that since $2hx^2-x^{2-\d}>hx^2$
 \bnn\label{eq_Llong}
 \left\{\inf_{0\le i\le h x^2}X_{\sigma_k+i}\ge x\sqrt{2h}\right\}
  \subseteq  \{(\sigma_k+i,X_{\sigma_k+i})\not\in L \mbox{ for all }0\le i\le  hx^2\}.
 \enn
By Lemma~\ref{lem_mart_clt}, the probability of the LHS
of~(\ref{eq_Llong}) is larger than
$$
1-\frac{4hB_1^2}{(1-\sqrt{2h})^2}=\frac 12.
$$

Similarly, when $(\eta_k,X_{\eta_k})\in R$ set
$y:=X_{\eta_k}>x^{\frac{2-\d}{2-2\d}}$. Since
 \bn
(y-x)^{2-\d}-x^{2-\d}>x^{\frac{(2-\d)^2}{2-2\d}}(1+o(1))-x^{2-\d}
 = x^{2+\frac{\d^2}{2-2\d}}(1+o(1))
 \gg x^2
 \en
 we have
 \bnn\label{eq_Rlong}
 \left\{\inf_{0\le i\le  x^2/(8B_1^2)}X_{\eta_k+i}\ge y-x\right\}
 \subseteq  \left\{(\eta_k+i,X_{\eta_k+i})\in M \mbox{ for all }0\le i\le
 \frac{x^2}{8B_1^2}\right\}.
 \enn
By Lemma~\ref{lem_mart_clt} the  probability of the LHS
of~(\ref{eq_Rlong}) is also less than $1/2$.  Therefore, since
$a<1/(8B_1^2)$, from (\ref{eq_sigma_large}) we obtain
 \bn
\P\left(\sigma_{k+1}-\sigma_k\ge a X_{\sigma_k}^2\right)>\frac 12.
 \en
On the other hand, provided $\sigma_k$ is large enough,
 \bn
a X_{\sigma_k}^2 =a \sigma_k^{\frac 2{2-\d}}(1+o(1))>3\sigma_k
 \en
yielding that for large $k$ for some $C_1>0$
 \bn
\sigma_{k}\ge C_1 4^{k/2}=C_1 2^k.
 \en
Consequently, the probability in (\ref{eq_p_l}) is bounded by
$$
 \frac {1+o(1)}{(C_1 2^k)^{\frac{\d}{2-\d}}},
$$
which is summable over $k$. By the Borel-Cantelli lemma only
finitely many events $\{(\eta_k,X_{\eta_k})\in L\}$ occur, or,
equivalently,  for large times $(t,X_t)\notin L$. This yields
transience.
 \Cox

\section{Recurrence}\label{sec_rec}
\begin{thm}\label{t2}
Consider a Markov process $X_t$, $t=0,1,2,\dots$ on $\R_+$ with
increments $D_t=X_{t}-X_{t-1}$, satisfying (H1) and (H2) for some
$a>0$. Suppose that on the event $\{X_t\ge a\}$  either
\begin{list} {(\roman{tally})} {\usecounter{tally}}
 \item for some $\rho<1/2$
 \bn
 \E(D_{t+1}\| \F_{t})\le \frac{\rho X_t}t,
 \en
 or
 \item for some $\rho>0$ and $(\a,\b)\in \Rec$
 \bn
 \E(D_{t+1}\| \F_{t})\le \frac{\rho X_t^{\a}}{t^\beta}.
 \en
\end{list}
Then $X_t$ is ``recurrent'' in the sense that  for any starting
point $X_0=x$ we have
$$
\P(\exists \, t\ge 0\text{ such that }X_t<a)=1.
$$
Hence also $\P(X_t<a\mbox{ infinitely often})=1$.
\end{thm}
\proof Consider  $ Y_t=X_t^2/t\ge 0$ and assume $X_t\ge a$. Then
 \bn
Y_{t+1}-Y_t=\frac{(X_t+D_{t+1})^2}{t+1}-\frac{X_t^2}{t}=
\frac{2tX_tD_{t+1}-X_t^2+t D_{t+1}^2}{t(t+1)}
 \en
whence
 \bnn\label{ekappa}
(t+1) \E(Y_{t+1}-Y_t\|\F_t)&=& \E(D_{t+1}^2\|\F_t)
+[2 X_t\E(D_{t+1}\|\F_t)-X_t^2/t] \nonumber \\
&\le& B_1^2 + \left(2\rho\kappa_t-1\right)X_t^2/t \le B_1^2
-\left(1-2\rho\kappa_t\right) Y_t
 \enn
where
$$
\kappa_t=\frac{X_t^{\a-1}}{t^{\b-1}}.
$$
Consider the following three cases:
\begin{itemize}
 \item[(i)] $\a=\b=1$, then $\kappa_t=1$;
 \item[(ii-a)] $\b>\a\ge 1$, then since $X_t\le B_1 t$,
$\kappa_t\le
 B_1^{\a}/t^{\b-\a}\to 0$ as $t\to\infty$;
 \item[(ii-b)] $\a<1$, $\b>(\a+1)/2$, then whenever
$Y_t=X_t^2/t\ge r$
 for some fixed positive constant $r$ we have
 \bn
  \kappa_t=\frac{1}{X_t^{1-\a}t^{\b-1}}\le  \frac{r^{(\a-1)/2}}{t^{(1-\a)/2}t^{\b-1}}
  =\frac{r^{(\a-1)/2}}{t^{\b-(\a+1)/2}}\to 0 \mbox{ as $t\to\infty$.}
 \en
\end{itemize}
(Note that (ii-a) and (ii-b) together cover the set $\Rec$.) Set
$r=B_1^2/(1-2\rho)>0$ in the first case, and set $r=2B_1^2$
otherwise. Then for $t$ sufficiently large from (\ref{ekappa}) we
have
 \bnn\label{eq_super}
 \E(Y_{t+1}-Y_t\|\F_t)   &\le& 0.
 \enn
Let $s\ge 0$, and set $\tau(s)=\inf\{t\ge s:\ Y_t\le r\}$.
Equation~(\ref{eq_super}) yields that $Y_{t\wedge \tau(s)}$ is a
supermartingale, therefore a.s.\ there is a limit
$Y_\infty=\lim_{t\to\infty}Y_{t\wedge \tau(s)}$. This implies that
either $\tau(s)<\infty$ for infinitely many $s\in\Z_+$, or that
there is a (random) $S$ such that $\tau(S)=\infty$. In both cases
we conclude that there is a possibly random value $Z$ such that
$X_t^2\le Z t$ for infinitely many times $t_k\in \Z_+$,
$k=1,2,\dots$.

First, suppose that $\a\ge 0$. Then for a fixed $t_k$ define a
process $X'_t=X_{t+t_k}$. Set $n=2 X_{t_k}$ and $\g=1/2$, and
observe that the process $X'_t$ satisfies the conditions of
Lemma~\ref{lem_snos} with some $c=c(2\b-\a,\rho,Z)>0$. Indeed,
when $X_t\le n$, the drift of $X_t$ is at most of order
$n^{\a}/t^{\b}\sim 1/n^{2\b-\a}\le 1/n$ since $2\b-\a\ge 1$ and
$\a\ge 0$. Hence, there is a constant $\nu>0$, independent of
$t_k$, such that
$$
\P(X_t,\ t\ge t_k,\ \text{reaches $[0,a]$ before
}[n,\infty)\|\F_{t_k})\ge \nu.
$$
Therefore, by the second Borel--Cantelli Lemma (Durrett, p.~240)
$\{X_t\le a\}$ for infinitely many~$t$'s.

\vskip 5mm

Now suppose that $\a<0$. Consider $W_t=X_t^{1-\nu}$ for some
$0<\nu<1$. Then
 \bnn\label{eq_W_sup}
  \E(W_{t+1}-W_t\|\F_t)&=& X_t^{1-\nu}  \E\left((1+D_{t+1}/X_t)^{1-\nu}-1\|\F_t\right)  \nonumber\\
 &=&(1-\nu)X_t^{1-\nu} \E\left(\frac{ D_{t+1}}{X_t}-\frac {\nu}{2}\frac{ D_{t+1}^2}{X_t^2}+O(X_t^{-3})\|\F_t\right)
   \nonumber \\
 &\le& (1-\nu) X_t^{-1-\nu} \left(\frac{ X_t^{1+\a}}{t^{\b}}-\frac {\nu B_2}{2}+O(X_t^{-3})\right)
 \enn
Let $n=n_k=\sqrt{Z t_k}$, so that $X_{t_k}\le n$. Since
$2\b>1+\a$, we can fix an $\zeta>1$ such that
$$
2\b>\zeta(1+\a).
$$
Consider the process $W_t$ for $t\in [t_k,\eta]$ where
$$
\eta=\eta_k:=\inf\{t\ge t_k:\ X_t\le a \mbox{ or } X_t\ge
n^{\zeta}\}.
$$
Then $\eta<\infty$ a.s. from the same argument as in part (i) of
Lemma~\ref{lem_snos}.

Moreover, for  $t\in [t_k,\eta]$
$$
\frac{ X_t^{1+\a}}{t^{\b}}\le
 \frac{  (n^\zeta)^{1+\a}}{(n^2/Z)^{\b}}= \frac {Z^\b}{n^{2\b-\zeta(1+\a)}}\to 0
 \mbox{ as } k\to\infty.
$$
since $t_k\to\infty$ and hence $n_k\to\infty$. Therefore the RHS
of (\ref{eq_W_sup}) is negative and thus $W_{t\wedge\eta}$ is a
supermartingale. Consequently, by the optional stopping theorem
$$
n^{1-\nu}\ge X_{t_k}^{1-\nu}=W_{t_k}\ge \E(W_\eta\|\F_{t_k})\ge
(1-p) (n^{\zeta})^{1-\nu}
$$
where $p=p_k=\P(X_{\eta}\le a\|\F_{t_k})$. This implies
 \bn
  p_k\ge 1-\frac 1{n_k^{(\zeta-1)(1-\nu)}}\to 1  \mbox{ as } k\to\infty.
 \en
finishing the proof of the theorem.
 \Cox

%
%
%

\section{Special cases}\label{sec_spec}
\subsection{Case $\a=\b\ge 0$}
Since we can always rescale the process $X_t$ by a positive
constant, in this section we assume that $B_1=1$.  Then, in turn,
it is also reasonable to restrict our attention only to the case
$\rho\le 1$, since if the jumps of $X_t$ can be indeed close to
$1$ with a positive probability, we might have $X\approx t$, and
the drift of order $\rho (X/t)^\b$ with $\rho>1$ would imply that
the drift is in fact larger than $1=B_1$ leading to a
contradiction, so the model would not be properly defined.

\begin{thm}[$\a=\b<1$]
Consider a Markov process $X_t$, $t=0,1,2,\dots$ on $\R_+$ with
increments $D_t=X_{t}-X_{t-1}$, satisfying (H1), (H2), and (H3)
for some $a>0$. Suppose that for some $\b<1$ and $\rho\in(0,1]$ on
the event $\{X_t\ge a\}$
 \bn
 \E(D_{t+1}\| \F_{t})\ge \rho \left(\frac{X_t}t\right)^\b,
 \en
Then $X_t$ is transient.
\end{thm}
\proof The proof is identical to the proof of Theorem~\ref{t1},
case (ii).
 \Cox

\begin{thm}[$\a=\b>1$]
Consider a Markov process $X_t$, $t=0,1,2,\dots$ on $\R_+$ with
increments $D_t=X_{t}-X_{t-1}$, satisfying (H1) and (H2) for some
$a>0$. Suppose that for some $\b>1$ and $\rho<1$  on the event
$\{X_t\ge a\}$ the process $X_t$ satisfies
 \bn
 \E(D_{t+1}\| \F_{t})\le \rho \left(\frac{X_t}t\right)^\b.
 \en
Then $X_t$ is recurrent.
\end{thm}
\proof Fix $\zeta\in(\rho,1)$ and consider $Y_t=X_t/t^{\zeta}$.
Then, calculating as before, we obtain
 \bn
 \E(Y_{t+1}-Y_t\|\F_t)\le
 \frac{X_t}{t(t+1)^\zeta}\left(\rho(X_t/t)^{\b-1}-\zeta\right).
 \en
Since $\limsup X_t/t\le B_1=1$ and $\rho<\zeta$, for large $t$
this is negative and hence $Y_t$ is a non-negative supermartingale
converging almost sure. On the other hand, $\zeta<1$, thus
implying
$$
 \lim_{t\to\infty} \frac{X_t}{t}=\lim_{t\to\infty} \frac {Y_t}{t^{1-\zeta}}=0 \text{  a.s.}
$$
and consequently since $\b>1$ for some sufficiently large $t$ we
have $(X_t/t)^{\b-1}<1/4$. Therefore, for large $t$,
 \bn
 \E(D_{t+1}\| \F_{t})\le \rho \left(\frac{X_t}t\right)^{\b-1}\times \frac{X_t}t \le \frac 14\frac{X_t}t,
 \en
and hence $X_t$ is recurrent by Theorem~\ref{t2}.
 \Cox

The following statement immediately follows from Theorems~\ref{t1}
and~\ref{t2}.
\begin{corollary}[$\a=\b=1$]
Suppose that $X_t$ is a process satisfying (H1), (H2), and (H3)
for some $a>0$.
\begin{list} {(\roman{tally})} {\usecounter{tally}}
 \item If for some $\rho<1/2$
 \bn
 \E(D_{t+1}\| \F_{t})\le  \frac{\rho X_t}{t} \ \mbox{ when $X_t\ge a$}
 \en
  then $X_t$ is recurrent.
 \item If for some $\rho>1/2$
 \bn
 \E(D_{t+1}\| \F_{t})\ge  \frac{\rho X_t}{t} \ \mbox{ when $X_t\ge a$}
 \en
then $X_t$ is transient.
\end{list}
\end{corollary}

\subsection{Case $\a\le 0$, $\b=0$}
In this case, the drift is of order $\rho/X_t^\nu$ where
$\nu=-\a\ge 0$. This is the situation resolved by
Lamperti~\cite{Lam1} and~\cite{Lam2}.
\begin{thm}[$\a=-1$, $\b=0$]
Suppose that $X_t$ is a process satisfying (H1) and (H2) for some
$a>0$. Then, when $X_t\ge a$,
\begin{list} {(\roman{tally})} {\usecounter{tally}}
 \item if for some $\rho\le 1/2$
 \bn
 \E(D_{t+1}\| \F_{t})\le \rho \frac{\E(D_{t+1}^2\|\F_t)}{X_t}
 \en
 then $X_t$ is recurrent;
 \item if for some $\rho> 1/2$
 \bn
 \E(D_{t+1}\| \F_{t})\ge \rho \frac{\E(D_{t+1}^2\|\F_t)}{X_t}
 \en
 then $X_t$ is transient.
\end{list}
\end{thm}
\begin{corollary}[$\a\in(-\infty,-1)\cup(-1,0)$,  $\b=0$]
Suppose that $X_t$ is a process satisfying (H1) and (H2) for some
$a>0$. Then, when $X_t\ge a$,
\begin{list} {(\roman{tally})} {\usecounter{tally}}
 \item if for some $\nu>1$
 \bn
 \E(D_{t+1}\| \F_{t})\le  \frac{\rho}{X_t^\nu}
 \en
 then $X_t$ is recurrent;
 \item if for some $\nu<1$
 \bn
 \E(D_{t+1}\| \F_{t})\ge  \frac{\rho}{X_t^\nu}
 \en
 then $X_t$ is transient.
\end{list}
\end{corollary}

\subsection{Case $2\b-\a=1$, $-1\le \a\le 1$: open problem}\label{subsec_open}
 Two cases $\a=\b=1$ and $\a=-1$,
$\b=0$ are already covered. It is also straightforward that when
$\a=0$, $\b=1/2$ by the law of iterated logarithm the process is
recurrent for any $\rho$.

Cases $-1<\a<0$, $\b=\frac 12 (\a+1)$ and $0<\a<1$, $\b=\frac 12
(\a+1)$: unfortunately, we cannot find a general sensible criteria
to separate recurrence and transience here, and leave this as an
{\bf open problem}.

\section{Application to urn models}\label{sec_urn}
Fix a constant $\sigma>0$. Consider a Friedman-type  urn process
$(W_n,B_n)$, with the following properties. We choose a white ball
with probability $W_n/(W_n+B_n)$ and a black ball with a
complementary probability; whenever we draw a white (black resp.)
ball, we add a random quantity $A$ of white (black resp.) balls
and $\sigma-A$ black (white resp.) balls; For simplicity, suppose
$0\le A\le \sigma$ a.s. A special case when $A$ is not random is
considered in Freedman (1964). Following his notations, let $\a=\E
A$, $\beta=\sigma-\a$, and $\rho=(\a-\b)/\sigma=(\a-\b)/(\a+\b)$.
Also assume that $\a>\b>0$.

It turns out that this urn can be coupled with a random walk
described above. Indeed, for $t=0,1,2,\dots$ set
$X_t=|W_t-B_t|/(\b-\a)\in \Z_+ \subset \R_+$. Without much loss of
generality assume that the process starts at time
$(W_0+B_0)/\sigma \in \Z$, then  $t=(W_t+B_t)/\sigma \in \Z$.

Consequently, once $X_t\ne 0$,
$$
\E(X_{t+1}-X_t\| \F_t)=
 \frac12\left(1+\frac{(\b-\a)X_t}{\sigma t}\right)(+1)
 +\frac12\left(1-\frac{(\b-\a)X_t}{\sigma t}\right)(-1) =
 \frac{\rho X_t}t.
$$
Corollary~3.3 in Friedman~(1965), states that when $\rho>1/2$,
$W_n-B_n=W_0-B_0$ (equivalently, $X_n=0$) occurs for finitely many
$n$ with a positive probability, and after the Corollary Friedman
says that he does not know whether this event has, in fact,
probability $1$. On the other hand, our Theorem~\ref{t2} answers
this question positively -- indeed, a.s.\ there will be finitely
many times when the difference between the number of white and
black balls in the urn equals a particular constant.

See Janson~\cite{Jan} and Pemantle~\cite{Review} for more on urn
models.

\end{document}